\documentclass[10pt]{amsart}
\usepackage{amsmath,amssymb,amsthm}

\begin{document}
\newtheorem{lem}{Lemma}[section]
\newtheorem{prop}{Proposition}[section]
\newtheorem{cor}{Corollary}[section]
\numberwithin{equation}{section}
\newtheorem{thm}{Theorem}[section]

\theoremstyle{remark}
\newtheorem{example}{Example}[section]
\newtheorem*{ack}{Acknowledgments}

\theoremstyle{definition}
\newtheorem{definition}{Definition}[section]

\theoremstyle{remark}
\newtheorem*{notation}{Notation}
\theoremstyle{remark}
\newtheorem{remark}{Remark}[section]

\newenvironment{Abstract}
{\begin{center}\textbf{\footnotesize{Abstract}}%
\end{center} \begin{quote}\begin{footnotesize}}
{\end{footnotesize}\end{quote}\bigskip}
\newenvironment{nome}

{\begin{center}\textbf{{}}%
\end{center} \begin{quote}\end{quote}\bigskip}

\newcommand{\triple}[1]{{|\!|\!|#1|\!|\!|}}

\newcommand{\xx}{\langle x\rangle}
\newcommand{\ep}{\varepsilon}
\newcommand{\al}{\alpha}
\newcommand{\be}{\beta}
\newcommand{\de}{\partial}
\newcommand{\la}{\lambda}
\newcommand{\La}{\Lambda}
\newcommand{\ga}{\gamma}
\newcommand{\del}{\delta}
\newcommand{\Del}{\Delta}
\newcommand{\sig}{\sigma}
\newcommand{\ome}{\omega}
\newcommand{\Ome}{\Omega}
\newcommand{\C}{{\mathbb C}}
\newcommand{\N}{{\mathbb N}}
\newcommand{\Z}{{\mathbb Z}}
\newcommand{\R}{{\mathbb R}}
\newcommand{\Rn}{{\mathbb R}^{n}}
\newcommand{\Rnu}{{\mathbb R}^{n+1}_{+}}
\newcommand{\Cn}{{\mathbb C}^{n}}
\newcommand{\spt}{\,\mathrm{supp}\,}
\newcommand{\SSS}{\mathcal{S}}
\newcommand{\F}{\mathcal{F}}
\newcommand{\xxi}{\langle\xi\rangle}
\newcommand{\eei}{\langle\eta\rangle}
\newcommand{\xei}{\langle\xi-\eta\rangle}
\newcommand{\yy}{\langle y\rangle}
\newcommand{\dint}{\int\!\!\int}
\newcommand{\hatp}{\widehat\psi}
\renewcommand{\Re}{\;\mathrm{Re}\;}
\renewcommand{\Im}{\;\mathrm{Im}\;}

\title[Decomposing $\dot H^s$-sequences in dispersive profiles]{The lack of compactness in the Sobolev-Strichartz inequalities}

\author{Luca Fanelli}
\address{Luca Fanelli: Universidad del Pais Vasco, Departamento de
Matem$\acute{\text{a}}$ticas, Apartado 644, 48080, Bilbao, Spain}
\email{luca.fanelli@ehu.es}

\author{Nicola Visciglia}
\address{Nicola Visciglia: Universit\`a di Pisa, Dipartimento di Matematica, Largo B.
Pontecorvo 5, 56100 Pisa, Italy}
\email{viscigli@dm.unipi.it}



\begin{abstract}
  We provide a general method to decompose any bounded sequence in $\dot H^s$ into linear dispersive profiles generated by an abstract propagator, with a rest which is small in the associated Strichartz norms. The argument is quite different from the one proposed by Bahouri-G\'erard and Keraani in the cases of the wave and Schr\"odinger equations, and is adaptable to a large class of propagators, including those which are matrix-valued.
\end{abstract}

\maketitle

\section{introduction}\label{sec:intro}
In the recent years, the research on nonlinear PDE's produced a relevant increment of strategies and techniques finalized to a complete understanding of some critical differential models. 
As a starting example, motivated by the interest on the Yamabe problem, some pioneer results were obtained by Aubin and Talenti in \cite{A,T}, giving answers to some natural questions related to the criticality of the Sobolev embedding $\dot H^s(\R^d)\subset L^{p(s)}(\R^d)$, with $p(s)=2d/(d-2s)$, and $0<s<d/2$. 
Some years later, a great and well celebrated contribution to the theory of critical elliptic PDE's was given by Pierre Louis Lions, who introduced the concentration-compactness method, which immediately turned out to be a standard tool (see \cite{L0, L1, L2}).
After the work by Lions, Solimini and G\'erard in \cite{G,S} independently, and with different proofs, were able to describe in a precise way the lack of compactness of the Sobolev embedding $\dot H^s(\R^d)\subset L^{p(s)}(\R^d)$ (and also the version for Lorentz spaces, in \cite{S}). Inspired to \cite{G}, Gallagher in \cite{Galla}, Bahouri and G\'erard in \cite{BG} and Keraani in \cite{K} proved analogous results related to the Sobolev-Strichartz estimates, respectively for the Navier-Stokes, the wave and the Schr\"odinger equation. As an example, we paste here the result proved by Keraani in \cite{K}: the following standard notations
\begin{equation*}
  L^p_tL^q_x:=L^p(\R;L^q(\R^d)),
  \qquad
  L^p_t\dot H^s_x:=L^p(\R;\dot H^s(\R^d)),
  \qquad
  L^r_{t,x}:=L^r_tL^r_x
\end{equation*}
will accompany the rest of the paper.
\begin{thm}[Keraani \cite{K}]\label{thm:keraani}
Let $(\varphi_n)_{n\geq0}$ be a bounded sequence in $\dot H^1(\R^3)$ and let $v_n(t,x):=e^{it\Delta}\varphi_n$.
Then
there exist a subsequence $(v'_n
)$ of $(v_n)$, a sequence $({\bf h^j})_{j\geq1}$, ${\bf h^j}=(h^j_n)_{n\geq0}$ for any $j\geq1$ of scales, a
sequence $({\bf z^j})_{j\geq1}=({\bf t^j,x^j})_{j\geq1}$, with ${\bf z^j}
=(t^j_n, x^j_n)_{n\geq0}$ for any $j\geq1$ of cores, and a sequence of functions $(U^j)_{j\geq1}$
in $\dot H^1(\R^3)$ such that: 
\begin{equation}\label{eq:kortog}
  \left |\frac{h^k_n}{h^j_n}\right|+
\left |\frac{h^j_n}{h^k_n}\right|+\left|\frac{t^j_n-t^k_n}{(h^j_n)^2}\right|+\left|\frac{x^j_n-x^k_n}{h^j_n}\right|
  \to+\infty,
\end{equation}
as $n\to\infty$, for any $j\neq k$;
\begin{equation}\label{eq:kdecomposition}
  v_n'(t,x)=\sum_{j=1}^l\frac{1}{\sqrt{h^j_n}}e^{i\left(\frac{t-t^j_n}{(h^j_n)^2}\right)\Delta}U^j
  \left(\frac{x-x^j_n}{h^j_n}\right)+w^l_n(t,x),
\end{equation}
for any $l\geq1$, with
\begin{equation}\label{eq:krest}
\limsup_{n\to\infty}\|w^l_n\|_{L^p_tL^q_x}\to0,
\end{equation}
as $l\to\infty$, for any (non-endpoint) $H^1$-admissible couple $(p,q)$ satisfying
\begin{equation*}
  \frac2p+\frac 3q=\frac 32-1,
  \qquad
  4< p\leq\infty;
\end{equation*}
\begin{equation}\label{eq:kenergy}
  \int|\nabla_x v_n'(0,x)|^2\,dx=
\sum_{j=1}^l  \int|\nabla U^j(x)|^2\,dx
  +\int|\nabla_xw_n^l(0,x)|^2\,dx+o(1),
  \qquad
  \text{as }n\to\infty.
\end{equation}
\end{thm}
Almost in the same years of \cite{BG}, \cite{K}, Kenig and Merle introduced in \cite{KM1, KM2} a new strategy to solve a large class of critical nonlinear Schr\"odinger and wave equations. The argument by Kenig and Merle is based on extrapolating, by contradiction, a single compactly behaving solution to the problem, which they call {\it critical element}, via concentration-compactness methods; then, the rigidity given by the algebra of the equation implies that such solution, with such compactness properties, cannot exist.
The basic tool in capturing the critical element is given by a nonlinear version of Theorem 
\ref{thm:keraani} (in the case of Schr\"odinger, and the analogous in \cite{BG} for wave), which is in fact a consequence of the same result and the scattering properties of the nonlinear flow. 
Since the Kenig-Merle proof turns out to be adaptable to a large class of nonlinear dispersive equations, a lot of results appeared in the very last years in the same spirit of Theorem 
\ref{thm:keraani}, for different propagators (see e.g. \cite{BV,BR, B, Galla, GKP, Koch, MV, R, RV}. Among the previous list, we mention the papers by Merle-Vega \cite{MV}, Begout-Vargas \cite{BV}, Rogers-Vargas \cite{RV} and recently Ramos \cite{R}, in which Strichartz estimates at the lowest scales are treated, and some refinements are needed, in the style of the one which has been proved by Moyua-Vargas-Vega in \cite{MVV}; in the cases of \cite{BG} and \cite{K} the inequality \eqref{eq:gerard} in the sequel, proved by G\'erard in 
\cite{G}, plays the analog role of the Strichartz refinement.

Therefore, it would be appreciable to have a general result, in the same style of Theorem 
\ref{thm:keraani}, which might hold for a large and unified class of dispersive propagators. On the other hand, as far as we can see, it is not clear if the strategy proposed in \cite{BG} and \cite{K} might be adaptable, in total generality, to many problems, as for example the case of dispersive systems.

In view of the above considerations, the aim of this paper is to provide a new proof, which is quite different from the one proposed in \cite{BG} and \cite{K}, and which works for a large amount of dispersive propagators, including among the others the matrix-valued cases.

We are now ready to prepare the setting of our main theorem.
In the following, we work with vector-valued functions $f=(f_1,...,f_N):\R^d\rightarrow \C^N$,
with the notation
\begin{equation}\label{eq:Hsobolev}
  \|f\|_{\dot H^s_x}^2=\sum_{j=1}^N \|f_j\|_{\dot H^s_x}^2.
\end{equation}
With the symbol $\mathcal L=\mathcal L(D)$ we denote an operator
\begin{equation*}
  \mathcal L(D)=\mathcal F^{-1}\left(\mathcal L(\xi)\mathcal F\right),
  \qquad
  \mathcal L(\xi)=(\mathcal L_{ij}(\xi))_{i,j=1\dots N}:\R^d\to\mathcal M_{N\times N}(\C),
\end{equation*}
where $\mathcal F$ is the standard Fourier transform, and the matrix $\mathcal L(\xi)$ is assumed to be hermitian; in the above setting, the dispersive character of the Cauchy problem
\begin{equation}\label{eq:eq}
  \begin{cases}
    i\partial_t u+\mathcal{L}(D) u=0
    \\
    u(0,x)=f(x)
  \end{cases}  
\end{equation}
just depends on the geometrical properties of the graph of $\mathcal L(\xi)$, as it is well known.
In addition, we make the following abstract assumptions:
\begin{itemize}
  \item[(H1)]
  there exists $0<s<\frac d2$ such that 
  the problem \eqref{eq:eq} is globally well-posed in $\dot H^s_x$, and 
the unique solution is given via the propagator $u(t,x)=e^{it\mathcal L(D)}f(x)$;
  \item[(H2)]  the flow $e^{it\mathcal L(D)}$ is unitary onto $\dot H^s_x$, i.e.
  \begin{equation*}
    \left\|e^{it\mathcal L(D)}f\right\|_{\dot H^s_x}=\|f\|_{\dot H^s_x},
    \qquad
    \forall t\in\R,
  \end{equation*}
  where $s$ is the same as in (H1), and $\|\cdot\|_{\dot H^s_x}$ is defined in \eqref{eq:Hsobolev}.
  \item[(H3)] the symbol $\mathcal L(\xi):\R^d\to\mathcal M_{N\times N}(\C)$ is $\alpha$-homogeneous, 
i.e., for all $\lambda>0$,
  \begin{equation*}
  \mathcal L(\lambda\xi)=\lambda^\alpha\mathcal L(\xi);
\end{equation*} 
\item[(H4)] there exist $2\leq p<q\leq\infty$ such that the following Strichartz estimate hold
\begin{equation*}
\left\|e^{it\mathcal L(D)} f\right\|_{L^p_tL^q_x}\leq C\| f\|_{\dot H^{s}_x},
\end{equation*}
with the same $s$ as in (H1) and some constant $C>0$.
\end{itemize}
By homogeneity, the couple $(p,q)$ in (H4) needs to satisfy the scaling condition
\begin{equation}\label{eq:admis}
\frac \alpha p+\frac dq=\frac d2-s,
\end{equation}
where $s$ is given by (H1) and $\alpha$ is the one in (H3).
Notice that, by the Sobolev embedding $\dot H^s_x\subset L^{\frac{2d}{d-2s}}_x$, 
for $0<s<d/2$, and the $\dot H^s_x$-preservation
\begin{equation*}
  \left\|e^{it\mathcal L(D)}f\right\|_{L^\infty_t\dot H^s_x}=\|f\|_{\dot H^s_x}
\end{equation*}
(assumption (H2) above),
we get
\begin{equation}\label{eq:1}
  \left\|e^{it\mathcal L(D)}f\right\|_{L^\infty_tL^{\frac{2d}{d-2s}}_x}\leq C\|f\|_{\dot H^s_x},
  \qquad
  0<s<\frac d2,
\end{equation}
for some constant $C>0$. Consequently, by interpolation with \eqref{eq:1}, 
an estimate as the one of assumption (H4) automatically holds for any $s$-admissible pair $(\tilde p, \tilde q)$, 
i.e. any $(\tilde p,\tilde q)$ satisfying \eqref{eq:admis}, with $\tilde p\geq p$.
In particular, we have
\begin{equation}\label{eq:est}
  \left\|e^{it\mathcal L(D)}f\right\|_{L^r_{t,x}}\leq C\|f\|_{\dot H^s_x},
  \qquad
  0<s<\frac d2,
  \qquad
  r=\frac{2(\alpha+d)}{d-2s}.
\end{equation}
There are several examples of operators $\mathcal L(D)$ satisfying the previous assumptions, including the cases of Schr\"odinger, non-elliptic Schr\"odinger, wave and Dirac propagators, as we will show later during the introduction.
We are now ready to state our main theorem.
\begin{thm}\label{thm:main}
  Let ${\bf u}=(u_n)_{n\geq 0}$ be a bounded sequence in $\dot H^s_x$ for $0<s<\frac d2$. There exist a subsequence $(u'_n
)$ of $(u_n)$, a sequence $({\bf h^j})_{j\geq1}$, ${\bf h^j}=(h^j_n)_{n\geq 0}$ of scales, for any 
$j\geq1$, a
sequence $({\bf z^j})_{j\geq1}=({\bf t^j,x^j})_{j\geq1}$ of cores, with ${\bf z^j}=(t^j_n, x^j_n)_{n\geq0}$ 
for any $j\geq1$, and a sequence of functions $(U^j)_{j\geq1}$ in $\dot H^s_x$ such that: 
\begin{equation}\label{eq:mainortog}
  \left |\frac{h^m_n}{h^j_n}\right|+\left |\frac{h^j_n}{h^m_n}\right|+
  \left|\frac{t^j_n-t^m_n}{(h^j_n)^\alpha}\right|+\left|\frac{x^j_n-x^m_n}{h^j_n}\right|
  \to+\infty,
\end{equation}
as $n\to\infty$, for any $j\neq m$;
\begin{equation}\label{eq:maindecomposition}
  u_n'(x)=\sum_{j=1}^J\frac{1}{(h^j_n)^{\frac d2-s}}e^{i\left(\frac{t^j_n}{(h^j_n)^\alpha}\right)
\mathcal L(D)}U^j
  \left(\frac{x-x^j_n}{h^j_n}\right)+R^J_n(x),
\end{equation}
where $\alpha$ is the one in (H3),
for any $J\geq1$, with
\begin{equation}\label{eq:mainrest}
\limsup_{n\to\infty}\|e^{it\mathcal L (D)}R^J_n(x)\|_{L^{\tilde p}_tL^{\tilde q}_x}\to0,
\end{equation}
as $J\to\infty$, for any couple $(\tilde p,\tilde q)$ satisfying the admissibility condition \eqref{eq:admis}, with $p<\tilde p<\infty$, and $p$ is the one given by (H4);\\
for any $J\geq 1$ we have \begin{equation}\label{eq:kenergyprime}
  \|u_n'(x)\|^2_{\dot H^s_x}=
\sum_{j=1}^J  \|U^j(x)\|_{\dot H^s_x}^2
  +\|R_n^J(x)\|^2_{\dot H^s_x}+o(1),
  \qquad
  \text{as }n\to\infty.
\end{equation}
\end{thm}
\begin{remark}\label{rem:k0}
  Notice that \eqref{eq:maindecomposition} is slightly different to 
  \eqref{eq:kdecomposition}; effectively, it is sufficient to act at both the sides of 
  \eqref{eq:maindecomposition} with the propagator $e^{it\mathcal L(D)}$, to obtain the analogous of \eqref{eq:kdecomposition}. In fact, we prefer to write \eqref{eq:maindecomposition} in this form, because it respect the 
  {\it stationary} character of our proof. As it will be clear in the sequel, the main difference with the argument in \cite{BG}, \cite{K} is that, at each step of the recurrence argument which permits to extract the final sequence $u'_n$, we work on fixed sequences of times; arguing in this way, all the construction can be performed exactly as in the stationary theorem by G\'erard in \cite{G}. This idea is suggested by the argument which has been introduced in \cite{FVV}, to prove the existence of maximizers for Sobolev-Strichartz inequalities.
\end{remark}
\begin{remark}\label{rem:k}
  Theorem \ref{thm:main} implies Theorem \ref{thm:keraani}, in the special cases 
$\mathcal L(D)=\Delta$, $s=1$, $d=3$, apart from \eqref{eq:mainrest}; indeed, the case $p=\infty$ is missing in \eqref{eq:mainrest}. We do not find possible to obtain the decay of the $L^\infty_tL^{p(s)}_x$ in total generality; on the other hand, 
it is possible to prove it case by case, using
each time the specific properties of $\mathcal L(\xi)$.
By the way, we stress that the decay of $L^r_{t,x}$-norm in \eqref{eq:mainrest}, when $r$ is the one in \eqref{eq:est}, is typically 
the only information which is needed in the nonlinear applications.
\end{remark}
We now pass to give some examples of applications of the main theorem to other types of propagators.
\begin{example}[Wave propagator]\label{ex:wave}
The Strichartz estimates for the wave propagator $e^{it|D|}$ 
(see \cite{GV}, \cite{KT}), in dimension $d\geq2$, are the following:
\begin{equation}\label{eq:striwave}
  \left\|e^{it|D|}f\right\|_{L^p_tL^q_x}\leq C\|f\|_{\dot H^{\frac 1p-\frac 1q+\frac 12}_x},
\end{equation}
under the admissibility condition 
\begin{equation}\label{eq:admiswave}
  \frac 2p+\frac{d-1}q=\frac{d-1}2,
  \qquad
  p\geq2.
  \qquad
  (p,q)\neq(2,\infty).
\end{equation}
The gap of derivatives $\frac 1p-\frac 1q+\frac 12\geq0$ is 
null only in the case of the energy estimate $(p,q)=(\infty,2)$. 
In particular, 
\begin{equation}\label{eq:striwave2}
  \left\|e^{it|D|}f\right\|_{L^{\frac{2(d+1)}{d-1}}_{t,x}}\leq C\|f\|_{\dot H^{\frac12}_x}
  \qquad
  d\geq2,
\end{equation}
which is in fact the original estimate proved by Strichartz in \cite{STR}.
More generally, by Sobolev embedding one also obtains that
\begin{equation}\label{eq:striwave20}
  \left\|e^{it|D|}f\right\|_{L^{\frac{2(d+1)}{d-1-2\sigma}}_{t,x}}\leq C\|f\|_{\dot H^{\frac12+\sigma}_x},
  \qquad
  0\leq\sigma
<\frac{d-1}2,
\qquad
  d\geq2.
\end{equation}
Theorem \ref{thm:main} applies in this case, for any dimension $d\geq2$, 
and $0<\sigma<\frac{d-1}2$; notice that the case $\sigma=0$ is not included, since in this case
assumption (H4) fails. The case $\sigma=0$ has been recently treated and solved by Ramos in \cite{R}. 
\end{example}
\begin{example}[Dirac propagator]\label{ex:dirac}
  In dimension $d=3$, the massless Dirac operator is given by
  \begin{equation*}
   \mathcal D:=\frac1i\sum_{j=1}^3\alpha_j\partial_j.
  \end{equation*}  
  Here $\alpha_1,\alpha_2,\alpha_3\in\mathcal M_{4\time4}(\C)$ are the so called {\it Dirac matrices}, 
which are $4\times4$-hermitian matrices, $\alpha_j^t=\overline{\alpha_j}$, $j=1,2,3$, 
with the explicit form
  \begin{equation*}
    \alpha_1=
    \left(
    \begin{array}{cccc}
    0 & 0 & 0 & 1
    \\
    0 & 0 & 1 & 0
    \\
    0 & 1 & 0 & 0
    \\
    1 & 0 & 0 & 0
    \end{array}
    \right),
    \ \ \ 
    \alpha_2=
    \left(
    \begin{array}{cccc}
    0 & 0 & 0 & -i
    \\
    0 & 0 & i & 0
    \\
    0 & -i & 0 & 0
    \\
    i & 0 & 0 & 0
    \end{array}
    \right),
    \ \ \ 
    \alpha_3=
    \left(
    \begin{array}{cccc}
    0 & 0 & 1 & 0
    \\
    0 & 0 & 0 & -1
    \\
    1 & 0 & 0 & 0
    \\
    0 & -1 & 0 & 0
    \end{array}
    \right);
  \end{equation*}
  equivalently, $\alpha_j=\left(\begin{array}{cc}0 & \sigma_j \\ \sigma_j & 0\end{array}\right)$, 
where $\sigma_j$ is the $j^{\text{th}}$ $2\times2$-Pauli matrix, $j=1,2,3$.
  
  Since $\mathcal D^2=-\Delta I_{4\times4}$, the Strichartz estimates for 
the massless Dirac operator are the same as for the 3D wave equation (see \cite{DF}):  
  \begin{equation}\label{eq:stridir}
    \left\|e^{it\mathcal D} f\right\|_{L^p_tL^q_x}\leq C\|
f\|_{\dot H^{\frac 1p-\frac 1q+\frac 12}_x}.
  \end{equation}
  Here, the admissibility condition reads as follows:
  \begin{equation}\label{eq:admisdir}
  \frac 2p+\frac{2}q=1,
  \qquad
  p>2.
\end{equation}
In particular we have
\begin{equation}\label{eq:stridir2}
  \left\|e^{it\mathcal D} f\right\|_{L^{\frac{8}{2-2\sigma}}_{t,x}}\leq C\| f\|_{\dot H^{\frac12+\sigma}_x},
  \qquad
  0\leq\sigma<1.
\end{equation}
Also in this case Theorem \ref{thm:main} applies; moreover, the statement 
also includes the cases of more general dispersive systems.  At our knowledge, this is not a known fact; 
indeed, it is unclear if it might 
be possible to prove the same using the arguments by Bahouri-G\'erard and Keraani in \cite{BG}, \cite{K}.
\end{example}
\begin{example}[Non-elliptic Schr\"odinger propagators]\label{ex:schro}
  Let us consider the Schr\"odinger operator 
  $L:=\sum_{j=1}^m\partial_j^2-\sum_{j=m+1}^d
  \partial_j^2$, for $d\geq2$ and $1\leq m<d$. The Strichartz estimates are the same as for the Schr\"odinger propagator, namely
  \begin{equation}\label{eq:strischro}
     \left\|e^{itL}f\right\|_{L^p_tL^q_x}\leq C\|f\|_{\dot H^s_x},
  \end{equation}
  with the admissibility condition
  \begin{equation}\label{eq:admisschro}
    \frac 2p+\frac dq=\frac d2-s,
    \qquad
    p\geq2,
    \qquad
    (p,q)\neq(2,\infty).
  \end{equation}
  In particular, one has
  \begin{equation}\label{eq:strischro2}
    \left\|e^{itL}f\right\|_{L^{\frac{2(d+2)}{d-2s}}_{t,x}}\leq C\|f\|_{\dot H^s_x},
  \end{equation}
  for any $0\leq s<\frac d2$. The only case in which at our knowledge has been treated is $d=2,s=0$. Indeed, the result by Rogers and Vargas in \cite{RV} contains all the ingredients which are necessary to prove the profile decomposition for bounded sequences in $L^2$, with respect to the propagator $e^{itL}$, in dimension $d=2$.
  It is a matter of fact that Theorem \ref{thm:main} applies for any $d\geq2$ and $0<s<\frac d2$, but we remark that it cannot include the case $s=0$. In addition, our argument is rather simple and does not involve any Fourier properties of the propagator, since it just use the fact that $s>0$.
\end{example}
The rest of the paper is devoted to the proof of Theorem \ref{thm:main}.

\section{Proof of Theorem \ref{thm:main}}\label{sec:proof}

Let us start with some preliminary definitions, introduced in \cite{G}.
\begin{definition}\label{def:h}
  Let ${\bf f}=(f_n)_{n\geq1}$ be a bounded sequence in $L^2(\R^d)$ and ${\bf h}=(h_n)_{n\geq1}, 
{\bf \widetilde h}=(\widetilde h_n)_{n\geq1}\subset\R$
  two scales. We say that:
  \begin{itemize}
  \item
  ${\bf f}$ is ${\bf h}$-oscillatory if
  \begin{equation}\label{eq:osc}
    \limsup_{n\to\infty}\left(\int_{h_n|\xi|\leq\frac1R}|\widehat f_n(\xi)|^2\,d\xi
    +(\int_{h_n|\xi|\geq R}|\widehat f_n(\xi)|^2\,d\xi\right)\to0
    \qquad
    \text{as }R\to\infty;
  \end{equation}
  \item
  ${\bf f}$ is ${\bf h}$-singular if, for every $b>a>0$, we have
  \begin{equation}\label{eq:sing}
    \lim_{n\to\infty}\int_{a\leq h_n|\xi|\leq b}|\widehat f_n(\xi)|^2\,d\xi
   =0
    \qquad
    \text{as }R\to\infty;
  \end{equation}
  \item
  ${\bf h}$ and ${\bf \widetilde h}$ are orthogonal if
  \begin{equation}\label{eq:defort}
    \lim_{n\to\infty}\left(\frac{h_n}{\widetilde h_n}+\frac{\widetilde h_n}{h_n}\right)=0.
  \end{equation}
  \end{itemize}
\end{definition}
The following proposition, proved in \cite{G}, permits to reduce the matters to prove Theorem 
\ref{thm:main} in the case of ${\bf 1}$-oscillating sequences.
\begin{prop}[G\'erard \cite{G}]\label{prop:ger}
  Let ${\bf f}=(f_n)_{n\geq0}$ be a bounded sequence in $L^2(\R^d)$. There exist a subsequence ${\bf f'}$ of ${\bf f}$, a family $({\bf h^j})_{j\geq1}$ of pairwise orthogonal scales and a family $({\bf g^j})_{j\geq1}$ of bounded sequences in $L^2(\R^d)$ such that:
  \begin{itemize}
  \item
  ${\bf g^j}$ is ${\bf h^j}$-oscillatory, for every $j$;
  \item
  for every $J\geq1$ and $x\in\R^d$,
  \begin{equation}\label{eq:decger}
    f_n'(x)=\sum_{j=1}^Jg_n^j(x)+R^J_n(x),
  \end{equation}
  where $(R^J_n)_{n\geq1}$ is ${\bf h^j}$-singular, for every $j=1,\dots,J$ and
  \begin{equation}\label{eq:besov}
    \limsup_{n\to\infty}\|R^J_n\|_{\dot B^0_{2,\infty}}\to0
    \qquad
    \text{as }J\to\infty;
  \end{equation}
  \item
  for every $J\geq 1$,
  \begin{equation}\label{eq:gerort}
    \|f_n'\|_{L^2}^2=\sum_{j=1}^J\|g_n^j\|_{L^2}^2+\|R^J_n\|_{L^2}^2+o(1),
    \qquad
    \text{as }n\to\infty.
  \end{equation}
  \end{itemize}
\end{prop}
We are now ready to prove the following result, which is the core of the proof of Theorem \ref{thm:main}.
\begin{prop}\label{prop:localized}
  Assume (H1)-(H2)-(H3)-(H4). Let ${\bf{u}}=(u_n)\subset\dot H^s_x$ 
be a ${\bf 1}$-oscillatory, bounded sequence in $\dot H^s_x$ with $0<s<\frac d2$. 
There exist a subsequence ${\bf{u'}}=(u_n')$ of ${\bf{u}}$, a family of cores 
$({\bf{z^j}})_{j\geq1}=(t^j_n, x^j_n)_{n\geq 0, j\geq1}\subset\R\times\R^d$, 
and a family of functions 
  $\left(U^j\right)_{j\geq1}$ in $\dot H^s_x$ such that:
  \begin{itemize}
  \item[(i)]
  for any $j\neq k$, we have
  \begin{equation}\label{eq:orthog}
  \left|t^j_n-t^k_n\right|
  +\left|x^j_n-x^k_n\right|\to\infty, 
  \qquad
  \text{as } n\to\infty;
  \end{equation}
  \item[(ii)]
  for all $J\geq1$ and $x\in\R^d$,
  \begin{equation}\label{eq:decomposition}
    u_n'(x)=\sum_{j=1}^Je^{it^j_n\mathcal L(D)}U^j(x-x^j_n)+R_n^J(x),
  \end{equation}
  with 
  \begin{equation}\label{eq:rest}
    \limsup_{n\to\infty}\left\|e^{it\mathcal L(D)}R_n^J(x)\right\|_{L^{\tilde p}_tL^{\tilde q}_x}\to0,
    \qquad
    \text{as }J\to\infty,
  \end{equation}
  for any $s$-admissible pair $(\tilde p,\tilde q)$, with $\tilde p>p$, and $p$ given by (H4).
  In addition, 
  \begin{equation}\label{eq:energy}
    \|u_n'(x)\|_{\dot H^s_x}^2=
    \sum_{j=1}^J\left\|U^j(x)\right\|_{\dot H^s_x}^2+\left\|R_n^J(x)\right\|_{\dot H^s_x}^2+o(1),
    \qquad
    \text{as }n\to\infty.
  \end{equation}
  \end{itemize}
\end{prop}

\begin{proof}

 Let us introduce the notation $S:=L^{\infty}_tL^{\frac{2d}{d-2s}}_x$, and recall \eqref{eq:1}.
  The proof is based on a construction by recurrence, which is quite different from the one used in \cite{BG}, \cite{K}.
  
Assume that $\liminf_{n\rightarrow \infty} \|e^{-it\mathcal L (D)}u_n(x)\|_S=0$,
  then \eqref{eq:decomposition} is satisfied provided that $u_n'$ is a subsequence
of $u_n$
such that
$\lim_{n\rightarrow \infty} \|e^{-it\mathcal L}u_n'(x)\|_S=0$, 
  $J=0$, $U^0\equiv0$ and $R^0_n\equiv u_n'$; 
 moreover in this case \eqref{eq:orthog} and \eqref{eq:energy} are trivially satisfied 
and the decay of the interpolated Strichartz norms \eqref{eq:rest} follows by
interpolation between 
the decay of the $S$-norm and the $L^p_tL^q_x$ a priori bound given by assumption (H4).
  Therefore we shall assume that
  \begin{equation*}
    \|e^{-it\mathcal L(D)}u_n(x)\|_S\geq2\delta_0,
  \end{equation*}
  for any $n\geq 0$ and some $\delta_0>0$. By the definition of $S$, there exists a 
sequence of times $(t^1_n)_{n\geq 0}\subset\R$ such that
  \begin{equation}\label{eq:prima}
    w_n^1(x):=e^{-it_n^1\mathcal L(D)}u_n(x),
    \qquad
    \|w_{n}^1(x)\|_{L_x^{\frac{2d}{d-2s}}}
    \geq\frac12\left\|e^{-it\mathcal L(D)}u_n(x)\right\|_S
    \geq\delta_0.
  \end{equation}
Arguibg as G\'erard in \cite{G}, we now denote by $\mathcal P({\bf
w^1})$ the set of all the possible 
weak limits in $\dot H^s_x$ of all the possible subsequences of 
$(w^1_n)$ with all their possible translations; moreover, let
 \begin{equation}\label{eq:gamma1}
   \gamma({\bf w^1}):=\sup\left\{\|\psi\|_{\dot H^s_x}:\psi\in\mathcal P({\bf
w^1})\right\}.
 \end{equation}
 As a consequence, there exist a sequence of centers $(x^1_n)_{n\geq 0}\subset\R^d$ 
and a subsequence of $u_n$ (that we still denote $u_n$) such that
 \begin{equation}\label{eq:seconda}
   e^{-it^1_n\mathcal L(D)}u_n(x+x^1_n)=
   w^1_n(x+x^1_n)
   \rightharpoonup
   U^1(x)
   \qquad
   \text{ weakly in }\dot H^s_x,
 \end{equation}
as $n\to\infty$, where
\begin{equation}\label{eq:ottava}
  \gamma({\bf w^1})\leq 2 \left\|U^1\right\|_{\dot H^s_x}.
\end{equation}
 Next we introduce $R_n^1(x)$ as follows
 \begin{equation}\label{eq:quinta}
    u_n(x)=e^{it^1_n\mathcal L(D)}U^1(x-x^1_n)+R^1_n(x),
 \end{equation}
and by \eqref{eq:seconda} we get 
 \begin{equation}\label{eq:sesta}
   e^{-it^1_n\mathcal L(D)}R^1_n(x+x^1_n)
   \rightharpoonup
   0
   \qquad
   \text{weakly in }\dot H^s_x,
 \end{equation}
 as $n\to\infty$.
Next notice that by combining \eqref{eq:seconda} with \eqref{eq:quinta}
and by recalling the definition of weak limit we deduce
\begin{equation}\label{eq:settima'}
    \left\|e^{-it^1_n\mathcal L(D)}R^1_n(x+x_n^1)\right\|_{\dot H^s_x}^2=
    \left\|e^{-it^1_n\mathcal L(D)}u_n(x+x^1_n) \right\|_{\dot H^s_x}^2
    -\left\|U^1(x)\right\|_{\dot H^s_x}^2
    +o(1),
\end{equation}
and by the $\dot H^s_x$-preservation (H2) implies
\begin{equation}\label{eq:settima}
    \left\|R^1_n(x)\right\|_{\dot H^s_x}^2=
    \left\|u_n(x)\right\|_{\dot H^s_x}^2
    -\left\|U^1(x)\right\|_{\dot H^s_x}^2
    +o(1).\end{equation}
Next assume that 
\begin{equation}\label{eq:liminf}
\liminf_{n\rightarrow \infty} \left\|e^{it\mathcal L(D)}R^1_n(x)\right\|_{S}=0
\end{equation} then \eqref{eq:decomposition},
\eqref{eq:energy} follow by \eqref{eq:quinta}, \eqref{eq:settima}, and \eqref{eq:rest}
follows by interpolation between the $S$-norm (that goes to zero on a suitable
subsequence due to
\eqref{eq:liminf}) 
with the Strichartz norm given 
by assumption (H4). Therefore, up to choose a subsequence, we can assume as before that
  \begin{equation*}
    \left\|e^{it\mathcal L(D)}R^1_n(x)\right\|_{S}\geq2\delta_1,
  \end{equation*}
  for any $n\geq0$ and some $\delta_1>0$. As a consequence, there exists a sequence
of times 
  $(t^2_n)_{n\geq 0}\subset\R$ such that
  \begin{equation}\label{eq:undicesima}
    w_n^2(x):=e^{-it_n^2\mathcal L(D)}R^1_n(x),
    \qquad
    \|w_{n}^2(x)\|_{L^{\frac{2d}{d-2s}}_x}
    \geq\frac12\left\|e^{-it\mathcal L(D)}R^1_n(x)\right\|_S
    \geq\delta_1.
  \end{equation}
  Define as above
   \begin{equation}\label{eq:gamma2}
   \gamma({\bf w^2}):=\sup\left\{\|\psi\|_{\dot H^s_x}:\psi\in\mathcal P({\bf
w^2})\right\}.
 \end{equation}
To this we associate a new sequence $(x^2_n)_{n\geq0}\subset\R^d$ and 
a subsequence of $R^1_n$, which we still call $R^1_n$, such that
 \begin{equation}\label{eq:nona}
   e^{-it^2_n\mathcal L(D)}R^1_n(x+x^2_n)=
   w^2_n(x+x^2_n)
   \rightharpoonup
   U^2(x)
   \qquad
   \text{weakly in }\dot H^s_x,
 \end{equation}
as $n\to\infty$; moreover we can assume that
\begin{equation}\label{eq:decima}
  \gamma({\bf w^2})\leq 2 \left\|U^2\right\|_{\dot H^s_x}.
\end{equation}
Next we introduce $R_n^2(x)$ as follows:
\begin{equation}\label{eq:quattordicesima'}
  R_n^1(x)=
  e^{it^2_n\mathcal L (D)}U^2(x-x^2_n)+
  R^2_n(x),
\end{equation}
By \eqref{eq:nona} we conclude that 
\begin{equation}\label{eq:quindicesima}
  e^{-it^2_n\mathcal L(D)}R_n^2(x+x^2_n)
   \rightharpoonup
   0
   \qquad
   \text{weakly in }\dot H^s,
\end{equation}
Moreover arguing as in \eqref{eq:settima}
we get
\begin{equation}\label{eq:sedicesima'}
  \left\|R^2_n(x)\right\|_{\dot H^s_x}^2=
  \left\|R^1_n(x)\right\|_{\dot H^s_x}^2
  -\left\|U^2(x)\right\|_{\dot H^s_x}^2
  +o(1).
\end{equation}
By combining \eqref{eq:quinta} and \eqref{eq:quattordicesima'} we obtain
\begin{equation}\label{eq:quattordicesima}
  u_n(x)=
  e^{it^1_n\mathcal L(D)}U^1(x-x^1_n)+
  e^{it^2_n\mathcal L(D)}U^2(x-x^2_n)+
  R^2_n(x),
\end{equation}
and
by combining \eqref{eq:sedicesima'} with \eqref{eq:settima} we get
\begin{equation}\label{eq:sedicesima}
  \left\|R^2_n(x)\right\|_{\dot H^s_x}^2=
  \left\|u_n(x)\right\|_{\dot H^s_x}^2
  -\left\|U^1(x)\right\|_{\dot H^s_x}^2
  -\left\|U^2(x)\right\|_{\dot H^s_x}^2
  +o(1).
\end{equation}
  The computations above describe an iterative procedure which at any 
  step $j=0,1,\dots$ permits to construct a (finite) family $U^1,\dots,U^j\in\dot
H^s_x$, a family of cores 
  $(t^1_n,x^1_n)_{n\geq1},\dots(t^j_n,x^j_n)_{n\geq0}\in\R\times\R^d$,
  and a sequence $R^j_n(x)\in\dot H^s_x$ such that (up to subsequence)  $u_n$ can be
written as
  \begin{equation}\label{eq:A}
  u_n(x)=e^{it_n^1\mathcal L(D)}U^1(x-x^1_n)
  +\cdots+e^{it_n^j\mathcal L(D)}U^j(x-x^j_n)
  +R_n^j(x),
  \end{equation} 
with the following extra properties:
\begin{equation}\label{eq:importante1} \|e^{-it^{j}_n\mathcal
L(D)}R^{j-1}_n(x)\|_{L^{\frac{2d}{d-2s}}_x}
    \geq\frac12\left\|e^{-it\mathcal L(D)}R^{j-1}_n(x)\right\|_S
    \geq\delta_{j-1}>0;
 \end{equation}  
\begin{equation}\label{eq:importante}
   e^{-it^{j}_n\mathcal L(D)}R^{j-1}_n(x+x^{j}_n)
   \rightharpoonup
   U^j(x)
   \qquad
   \text{ weakly in }\dot H^s_x;
\end{equation}
\begin{equation}\label{eq:importante'}
 \gamma ({\bf w^j})\leq 2 \|U_j\|_{\dot H^s_x}
\end{equation}
(where the sequence $(w^j_n)$ is defined by $e^{-it^{j}_n\mathcal L(D)}R^{j-1}_n(x))$
and $\gamma({\bf w^j})$ is defined according to \eqref{eq:gamma1});
\begin{equation}\label{eq:B}
  \left\|R^j_n(x)\right\|_{\dot H^s_x}^2=
  \left\|u_n(x)\right\|_{\dot H^s_x}^2
  -\left\|U^1(x)\right\|_{\dot H^s_x}^2
  -\cdots
  -\left\|U^j(x)\right\|_{\dot H^s_x}^2
  +o(1);
\end{equation} 
\begin{equation}\label{eq:0}
e^{-it^j_n\mathcal L(D)}R^j_n(x+x^j_n)\rightharpoonup0,
\hbox{ weakly in } \dot H^s_x.
\end{equation}
Notice that \eqref{eq:A} and \eqref{eq:B} prove \eqref{eq:decomposition} and
\eqref{eq:energy}. 
Our goal is now to prove \eqref{eq:rest}; it is sufficient to prove that, 
for any $\epsilon>0$, there exists $J=J(\epsilon)\in\mathbb N$ such that, 
for any $n\in\mathbb N$, and for any $j\geq J(\epsilon)$ we have
\begin{equation}\label{eq:claim}
\limsup_{n\rightarrow \infty}  \left\|e^{it\mathcal L(D)}R^j_n(x)\right\|_{S}<\epsilon.
\end{equation}
Fix $\epsilon>0$; first observe that, by \eqref{eq:B}
the sum $\sum_{j\geq 1}\|U^j\|_{\dot H^s_x}^2$ has to converge, 
then there exists $J=J(\epsilon)$ such that, for any $j\geq J(\epsilon)$,
\begin{equation}\label{eq:quasi}
\gamma({\bf w^{j}})\leq 2 
\left\|U^{j}\right\|_{\dot H^s_x}<2\epsilon
\end{equation}
where we have used \eqref{eq:importante'}.
In order to conclude \eqref{eq:claim} it is sufficient, by 
\eqref{eq:importante1}, to prove that
\begin{equation}\label{eq:quasi2}
\|e^{-it^{j}_n\mathcal L(D)}R^{j-1}_n(x)\|_{L^{\frac{2d}{d-2s}}_x}:=
  \|w^{j}_n(x)\|_{L^{\frac{2d}{d-2s}}_x}<C\epsilon,
\end{equation}
for any $j\geq J(\epsilon)$, and some constant $C>0$.
This is an immediate consequence of the inequality
\begin{equation}\label{eq:gerard}
  \limsup_{n\to\infty}\|w_n^{j}\|_{L^{p(s)}_x}\leq
  C\limsup_{n\to\infty}
  \left\|w_{n}^{j}\right\|_{\dot H^s_x}^{\frac2{p(s)}}\gamma({\bf
w^{j}})^{1-\frac2{p(s)}},
\end{equation}
with $p(s)=2d/(d-2s)$. The previous estimate has been proved by G\'erard (see \cite{G}, estimate (4.19)).

In order to complete the proof, we need to show the orthogonality of the cores
\eqref{eq:orthog}. 
Let us first prove it in the case $k=j+1$.
Notice that by \eqref{eq:importante} 
we have
\begin{equation}\label{eq:importante''}
   e^{-it^{j+1}_n\mathcal L(D)}R^{j}_n(x+x^{j+1}_n)
   \rightharpoonup
   U^{j+1}(x)
   \qquad
   \text{ weakly in }\dot H^s_x
\end{equation}
which is equivalen to
\begin{equation*}
e^{-i\left(\left(t^{j+1}_n-t^j_n\right)+t^j_n\right)\mathcal
L(D)}R^j_n\left(x+\left(x^{j+1}_n-x_n^j\right)+x^j_n\right)
  \rightharpoonup U^{j+1}(x).
\end{equation*}
Next assume by the absurd that the cores $({\bf z^j})$ and $({\bf z^{j+1}})$
do not satisfy \eqref{eq:orthog}, then up to subsequence we can assume
$t^{j+1}_n-t^j_n\rightarrow \bar t$
and $x^{j+1}_n-x_n^j\rightarrow \bar x$,
which in turn implies
\begin{equation*}
e^{-i(\bar t +t^j_n)\mathcal L(D)}R^j_n (x+\bar x +x^j_n )
  \rightharpoonup U^{j+1}(x).
\end{equation*}
Notice that this last fact is in contradiction with \eqref{eq:0}.\\
Next we assume that there exist a couple $(k, j)$ such that $k<j-1$ 
and for which the orthogonality 
\eqref{eq:orthog} for the cores $(\bf {z^j}), ({\bf z^k})$ is false (the case $k=j-1$ has been treated above).
Moreover we can suppose that 
the orthogonality relation is satisfied for
the cores $({\bf z^{k+r}})$ and $({\bf z^j})$ for any $r=1,..., j-k-1$. 
In fact it is sufficient to choose $k$ as 
$$\sup \{h<j-1|({\bf z^h}) \hbox{ is not orthogonal to } (\bf{z^j})\}.$$
Next notice that
$$R^k_n(x)= \sum_{h=k+1}^{j-1} e^{i t^h_n\mathcal L(D)} U_h(x-x^h_n)+R^{j-1}_n(x)$$
(to prove this fact apply \eqref{eq:A} twice: first up to the reminder $R^{k}_n$
and after up to the reminder $R^{j-1}_n$, and subtract the two identities). 
As a consequence we get
\begin{equation}\label{eq:identita} e^{-i t^j_n\mathcal L(D)} R_n^k(x+x_n^j)
\end{equation}$$=
\sum_{h=k+1}^{j-1} e^{i (t^h_n-t^j_n)\mathcal L(D)} U_h(x+x_n^j -x^h_n)+
e^{-i t^j_n \mathcal L(D)} R^{j-1}_n(x+ x_n^j).$$
Next notice that by the orthogonality 
of $({\bf z^{k+r}})$ and $({\bf z^j})$ for $r=1,...,j-k-1$ we get 
\begin{equation}\label{eq:identita2}
e^{i (t^h_n-t^j_n)\mathcal L} U_h(x+x_n^j -x^h_n) \rightharpoonup 0
 \end{equation}
 for every $h=k+1,...,j-1$ (here we use Lemma \ref{lem:fina}).
On the other hand we have the following identity
$$e^{-i t^j_n\mathcal L(D)} R_n^k(x+x_n^j)
= e^{-i (t^j_n-t^k_n+t^k_n)\mathcal L(D)} R_n^k(x+x_n^j-x_n^k+x_n^k)$$
and since we are assuming that $({\bf z^j})$ and $({\bf z^k})$ are not orthogonal, then
by compactness we can assume that $x_n^j-x_n^k\rightarrow \bar x$ and 
$t^j_n-t^k_n\rightarrow \bar t$. In particular 
we get $$e^{-i t^j_n\mathcal L(D)} R_n^k(x+x_n^j)
- e^{-i (\bar t +t^k_n)\mathcal L(D)} R_n^k(x+ \bar x +x_n^k)\rightarrow 0$$
and since by \eqref{eq:0} we have 
$e^{-i t^k_n\mathcal L(D)} R_n^k(x+x_n^k)\rightharpoonup 0$,
then necessarily also
the l.h.s. in \eqref{eq:identita} converges weakly to zero.
By combining this fact with \eqref{eq:identita2} we deduce that
$$e^{-i t^j_n \mathcal L(D)} R^{j-1}_n(x+ x_n^j)\rightharpoonup 0$$
and it is in contradiction with \eqref{eq:importante}.

\end{proof}

Having in mind Proposition \ref{prop:ger}, to complete the proof of Theorem \ref{thm:main} it is now sufficient to follow exactly the arguments given by Keraani in \cite{K}. One should only be careful at the moment of proving \eqref{eq:mainrest}; indeed, notice that in Lemma 2.7 in \cite{K} it is used the fact that $r=10$ is an integer number. On the other hand, the reader should easily notice that this is not a relevant fact, and the proof can be easily performed in the general case in which $r$ is given by \eqref{eq:est}. Once the decay of the $L^r_{t,x}$ of the rest is proved, the decay of the norms in \eqref{eq:mainrest} follows by interpolation. We omit here further details.

\appendix
\section{}
\label{sec:appendix}

We devote this small appendix to prove a general result, Lemma \ref{lem:fina}, 
which has been implicitly used during the proof of Proposition \ref{prop:localized}. Let us start with the following proposition.

\begin{prop}\label{prop:1}
Assume that 
\begin{equation}\label{strip2}
\left\|e^{it\mathcal L(D)} f\right\|_{L^p_tL^q_x}\leq C\| f\|_{\dot H^{s}_x},
\end{equation}
for some $p,q\geq1$, $p\neq \infty$ and some $C>0$.
Then $\left\|e^{it\mathcal L(D)} f\right\|_{L^q_x}
\rightarrow 0$, as $t\rightarrow \infty$, for any $f\in C^\infty_0(\R^d)$.
\end{prop}
\begin{proof}
Let $M>>1$ such that $H^M_x\subset L^q_x$, and let $f\in H^M_x$.
Then, by continuity in time,
for every $\epsilon>0$ there exists $\bar t=\bar t(\epsilon, f)>0$ such that
\begin{equation}\label{continuity}
\|e^{it\mathcal L(D)}f - f\|_{H^M_x}\leq \epsilon, 
\qquad
 \forall |t|<\bar t.
\end{equation}
Now assume by the absurd that
for a sequence $t_n\rightarrow \infty$ we have
\begin{equation}\label{absurd}\inf_n \left\|e^{it_n\mathcal L(D)} f\right\|_{L^q_x}=\delta>0.\end{equation}
As a consequence, by combining the Sobolev embedding
$H^M_x\subset L^q_x$ with \eqref{continuity} 
and the fact that  $e^{it\mathcal L(D)}$
is an isometry on $H^M_x$, we have
\begin{align*}
\|e^{i(t_n+h)\mathcal L(D)} f - e^{it_n\mathcal L(D)} f\|_{L^q_x}
&
\leq C \|e^{i(t_n+h)\mathcal L(D)} f - e^{it_n\mathcal L(D)} f\|_{H^M_x}
\\
& = C \|e^{ih\mathcal L(D)} f - f\|_{H^M_x}
\leq \frac \delta 2,
\end{align*}
provided that $|h|\leq h(\delta, f)$.
Therefore we deduce by
\eqref{absurd} that 
\begin{equation*}
\|e^{i(t_n+h)\mathcal L(D)} f\|_{L^q_x}\geq  \frac{\delta}2,
\end{equation*}
for any $n\in \N$ and $|h|\leq \bar h$.
The last estimate is in contradiction with \eqref{strip2}
since it does not allow global summability in time.
\end{proof}
We can now prove the main result of the appendix.
\begin{lem}\label{lem:fina}
Assume that
\begin{equation}\label{strip}
\left\|e^{it\mathcal L(D)} f\right\|_{L^p_tL^q_x}\leq C\| f\|_{\dot H^{s}_x},
\end{equation}
for some $p,q\geq1$ and some $C>0$.
Let $f\in \dot H^s_x$ and $\max\{|t_n|, |x_n|\}\rightarrow \infty$
then $e^{it_n\mathcal L(D)} f(x+x_n)\rightharpoonup 0$ in $\dot H^s_x$.
\end{lem}
\begin{proof}
We need to consider two cases. The first possibility is that $t_n$ is bounded; then 
necessarily $x_n$ goes to $\infty$ and it is easy to conclude.
In the case $t_n\rightarrow \infty$, the conclusion is now simple, by combining  a density argument
with Proposition \ref{prop:1}.
\end{proof}


\end{document}